\documentclass[10pt]{amsart}
\usepackage{amssymb}
\usepackage{epsfig}
\usepackage{amsfonts}
\usepackage{amscd}
\usepackage{amsmath}

\newcommand{\wo}{\overline{w}_0}
\newcommand{\g}{\Gamma(u,w;f)}

\newcommand{\mbotreix}{\overline{\mathcal{M}}_{0,3}(X,d)}

\newcommand{\z}{\mathbb{Z}}

\newcommand{\cx}{\mathbb{C}}
\newtheorem{defn}{Definition}[section]
\newtheorem{prop}{Proposition}[section]
\newtheorem{thm}[prop]{Theorem}
\newtheorem{cor}[prop]{Corollary}
\newtheorem{lemma}[prop]{Lemma}
\newtheorem{ithm}{Theorem}
\newcommand{\eqqhx}{{QH}^{*}_{T}(X)}
\newcommand{\pnop}{\Phi^+ \smallsetminus \Phi^+_P}
\newcommand{\dnop}{\Delta \smallsetminus \Delta_P}
\numberwithin{equation}{section}

\begin{document}

\title{Equivariant quantum cohomology of homogeneous spaces}
\author{Leonardo Constantin Mihalcea } \date{\today}
\begin{abstract} We prove a Chevalley formula for the equivariant
quantum multiplication of two Schubert classes in the homogeneous
variety $X=G/P$. As in the case when $X$ is a Grassmannian
(\cite{Mi1}), this formula implies an algorithm to compute the
structure constants of the equivariant quantum cohomology algebra
of $X$.
\end{abstract}

\address{University of Michigan, Dept. of Mathematics, East Hall,
525 E. University, Ann Arbor MI 48109-1109}
\email{lmihalce@umich.edu} \maketitle

\section{Introduction}

The aim of this paper is to prove a multiplication formula between
two Schubert classes in the equivariant quantum (EQ) cohomology of
a homogeneous variety $X=G/P$, where $G$ is a complex, semisimple,
connected linear algebraic group and $P$ a parabolic subgroup,
when one of these classes corresponds to a divisor. Despite the
fact that the (quantum, equivariant or EQ) cohomology of $X$ is
not in general generated by divisor classes, such a multiplication
will be sufficient for deriving an algorithm to compute all the
structure constants of the EQ cohomology algebra of $X$. This
algorithm implies in particular new algorithms to compute the
structure constants of the equivariant and, especially, quantum
cohomology of $X$, cohomologies which have been extensively
studied recently (see e.g. \cite{W,KM,AS,Be,FP,FGP,C,Ch,Bu3,BKT1}
and references therein for quantum cohomology and e.g.
\cite{A,Br1,GKM,EG,Ku} for equivariant cohomology).

Fix $T \subset B \subset P$, a Borel subgroup $B$ contained in the
parabolic subgroup $P$ together with its maximal torus $T$. The EQ
cohomology of $X$ is a deformation of both the equivariant and
quantum cohomology of $X$. It is a graded algebra over
$\Lambda[q]$, where $\Lambda$ denotes the $T-$equivariant
cohomology of the point and $q=(q_\beta)$ is a sequence of
indeterminates indexed by the simple reflections in the Weyl group
$W$ which are not in the Weyl group $W_P$ of $P$. It is well-known
that $\Lambda$ can be identified with the polynomial ring in the
negative\begin{footnote}{We have chosen the negative simple roots
instead of the positive ones to be the generators of $\Lambda$ for
positivity reasons: the structure constants in EQ cohomology are
nonnegative combinations of monomials in the {\em negative} simple
roots (\cite{Mi2}).}\end{footnote} simple roots of $G$, regarded
as characters of $T$. The EQ cohomology of $X$ has a
$\Lambda[q]-$basis consisting of Schubert classes $\sigma(u)$
indexed by the minimal length representatives for the cosets of
$W/W_P$. The multiplication of two Schubert classes is given by
\[ \sigma(u) \circ \sigma(v) = \sum_d \sum_w q^d c_{u,v}^{w,d}
\sigma(w) \] where $c_{u,v}^{w,d}$ are the (3-point, genus 0)
equivariant Gromov-Witten invariants introduced by Givental and
Kim (\cite{GK}). Here $d=(d_\beta)$ is a multidegree in $X$ (i.e.
a sequence of nonnegative integers), with the same number of
components as the indeterminate $q$, and $q^d$ stands for the
monomial $\prod q_\beta^{d_\beta}$ (all of these are detailed in
sections \ref{schubert}-\ref{eqqcohs} below). We will refer to the
coefficients $c_{u,v}^{w,d}$ as the equivariant quantum
Littlewood-Richardson coefficients (EQLR). They are homogeneous
polynomials in $\Lambda$, of degree \[c(u)+c(v)-c(w) - \sum
d_\beta \cdot (\deg q_\beta)\] where $c(u)$ is equal to the
complex degree of the cohomology class $\sigma(u)$ (see \S
\ref{schubert} below). The coefficient $c_{u,v}^{w,d}$ is equal to
the ordinary, non-equivariant, 3-point Gromov-Witten (GW)
invariant $c_{u,v}^{w,d}$ (which counts rational curves in $X$
subject to certain conditions) if its polynomial degree is equal
to zero, and to the structure constant $c_{u,v}^w$ of the
equivariant cohomology if its (multi)degree $d$ is equal to zero.
An important property of the EQLR coefficients is a certain
positivity, which generalizes the positivity enjoyed by the
structure constants of the equivariant cohomology (cf. \S
\ref{eqcohs}). In its original (equivariant) formulation it was
conjectured by Peterson and proved by Graham \cite{Gr}. The author
has proved the quantum generalization in \cite{Mi2}.

Givental and Kim \cite{GK,Kim1,Kim3} have used the EQ cohomology
successfully to prove properties about the (non-equivariant)
quantum cohomology, then Givental \cite{G} have used related ideas
to study Mirror Symmetry for projective complete intersections. In
fact, our initial motivation to carry out this work was to gain a
better understanding of the quantum cohomology of the homogeneous
spaces. It turns out however, that the EQ cohomology possesses
nice properties in its own right, therefore deserving a closer
study.

We will make distinction between the coefficients $c_{u,v}^{w,d}$
that are either quantum or equivariant ones, for which some
properties are known, and those that are not, which we call {\bf
mixed}. Thus the mixed coefficients are those for which neither
the polynomial degree nor the degree $d$ is equal to zero. In
general, small examples, see e.g. the EQ multiplication table for
the Grassmannian of $2-$planes in $\cx^4$ in \S 8, \cite{Mi1},
show that most of the mixed EQLR coefficients are nonzero. However

\begin{ithm} The mixed coefficients appearing in the EQ multiplication
with a Schubert class corresponding to a divisor (i.e. in a
Chevalley formula) are equal to zero. \end{ithm}

The precise multiplication formula is given in Thm. \ref{eqqc}.
The equivariant coefficients appearing in this formula have been
computed by Kostant and Kumar (\cite{KK}, see also
\cite{Bi})\begin{footnote}{A more general multiplication in
equivariant cohomology of the type A flag manifold, with cycle
permutations, has been obtained in \cite{R}.}\end{footnote} (in
fact, the coefficients computed there were those coming from a
basis of a restriction of a dual of the nil-Hecke algebra. The
relationship of this algebra with the equivariant cohomology
algebra of the flag varieties for the Kac-Moody groups, which
generalizes our situation, was established later by Arabia
\cite{A}, see also \cite{Ku}, Ch. 11). The formula for the quantum
coefficients (i.e. the quantum Chevalley rule) has been
conjectured by Peterson \cite{P} and proved by Fulton and Woodward
\cite{FW}.

As in the case when $X$ is a Grassmannian, studied by the author
in \cite{Mi1}, the theorem follows from a more general vanishing
property of the EQLR coefficients. However, Buch's notions of span
and kernel of a rational curve (see \cite{Bu1}), which were the
key tools used in the proof from \cite{Mi1}, are no longer
available for an arbitrary homogeneous space $G/P$. Nevertheless,
a slightly weaker vanishing (see Lemma \ref{main} below), which is
enough for our purposes, can be proved using a counting argument
using the Kleiman transversality Theorem. The analogy with the
Grassmannian goes further, as it turns out that the properties
characterizing the EQLR coefficients from \cite{Mi1} extend to
this more general context. More precisely
\begin{ithm}\label{ialg} The EQLR coefficients are uniquely determined
by the following:
\begin{enumerate}
\item[(a)](homogeneity) $c_{u,v}^{w,d}$ is a homogeneous {\em
rational} functions of degree \[c(u)+c(v)-c(w)- \sum d_\beta \cdot
(\deg q_\beta).\]
 \item[(b)](multiplication by unit)
$$ c_{id,w}^{w,d}= \left \{
\begin{array}{ll} 1 & \textrm{if } d=0\\
0 & \textrm { } \textrm{otherwise}
\end{array} \right.$$
where $0$ denotes the degree with all components equal to zero and
$\sigma(id)$ is the unit element for the EQ multiplication.

\item[(c)] (commutativity) $c_{u,v}^{w,d}=c_{v,u}^{w,d}$ for any
$u,v,w$.

\item[(d)] The equivariant quantum Chevalley coefficients, for
which a formula is given in (\ref{eqqcf}) below.

\item[(e)] A recurrence formula (see Cor. \ref{rec}) implied by
the associativity equation
\begin{equation}\label{assoc}
\bigl(\sigma(s(\beta)) \circ \sigma(u)\bigr) \circ \sigma(v) =
\sigma(s(\beta)) \circ \bigl(\sigma(u) \circ \sigma(v)\bigr)
\end{equation} where $\sigma(s(\beta))$ is a divisor Schubert
class (see \S \ref{schubert} below).

\end{enumerate}
\end{ithm}
The statement of the theorem is in fact slightly stronger (see
Thm. \ref{alg} below), and it is proved by exhibiting an effective
algorithm to compute the EQLR coefficients. Consequences of this
theorem include new algorithms for the computation of equivariant
and quantum coefficients. In the equivariant case, different
algorithms, and closed formulae for some of the coefficients have
been obtained in \cite{KK,Bi,K}. The only situation when closed
formulae for {\it all} of the coefficients are known is when $X$
is a Grassmannian (of type A). In this situation, a positive
formula (in the sense of \cite{Gr} or \cite{Mi2}, see \S
\ref{eqcohs} below) have been obtained by Knutson and Tao in
\cite{KT1}. Their formula is proved using a recursive expression,
which is in fact the recursive formula from (e) of Thm. \ref{ialg}
above, considered in equivariant cohomology.

Algorithms to compute the GW invariants in type A have been
obtained e.g. in \cite{Be,BCF,Po,Co} for Grassmannians, in
\cite{C,FGP} for complete flag manifolds and in \cite{C2,Bu3} for
partial flag manifolds. Some limited cases in other types have
been obtained in \cite{KTa,BKT1}. Since the quantum cohomology is
not functorial, the study of the GW invariants on an arbitrary
homogeneous variety $G/P$ is more difficult, and it has to be done
case by case. A very useful tool for this situation is the
Peterson comparison formula \cite{P}, proved by Woodward
\cite{Wo}, which shows that any GW invariant on $G/P$ is equal to
one on $G/B$, of possible different degree. The study of the
quantum cohomology of $G/B$ is simpler. In this case the
cohomology is generated by divisor classes, and a presentation, in
the more general context of EQ cohomology, has been obtained by
Kim \cite{Kim3} (see also \cite{Ma2}), in terms of the motions of
the Toda lattice. However, in order to compute the GW invariants
one needs to have polynomial representatives for the quantum
Schubert classes expressed in terms of these generators (i.e.
quantum Giambelli formulae) and to be able to multiply the
Schubert polynomials with the given generators (a quantum
Chevalley formula\begin{footnote}{In fact, S. Fomin pointed out
that given the quantum Giambelli formula, the quantum Chevalley is
not necessary for computational purposes; nevertheless, the
quantum Chevalley rule decreases dramatically the number of the
computations involved.}\end{footnote}).

In type A, a quantum Giambelli and Chevalley formula has been
obtained by Fomin-Gelfand-Postnikov \cite{FGP}. In other types, an
algorithm for the quantum Schubert polynomials has been obtained
by Mare in \cite{Ma1}, using ideas from \cite{FGP}. A
differential-geometric approach to the quantum cohomology of
$G/B$, initiated by Guest \cite{G}, was used recently to produce
algorithms for the quantum Schubert polynomials and to recover the
Chevalley formula (see also \cite{AG, Ma3}).

Our algorithm to compute the usual GW-invariants on $G/P$, implied
by the algorithm computing the EQLR coefficients, is conceptually
simpler. It needs the equvariant version of the GW invariants, but
just those invariants associated to $G/P$, hence avoiding
Peterson's comparison formula. It does not require the knowledge
of a presentation or EQ Giambelli formula and it uses only a
multiplication formula with divisor Schubert classes. Note that,
unless $P=B$, the divisors do not generate the EQ (or quantum)
cohomology algebras. Thus, in a certain sense,
the EQ cohomology behaves better than the quantum one.\\

{\it Acknowledgements:} I would like to thank my advisor, Prof. W.
Fulton for his encouragement and support and for his comments on
earlier versions of this paper. I would also like to thank S.
Fomin, B. Ion, L. Mare and M. Mizerski for some useful
conversations and helpful e-mail. This is part of my thesis.

\section{Cohomology of $G/P$}\label{schubert}

Throughout this paper $G$ denotes a complex, connected,
semisimple, linear algebraic group, $B$ a (fixed) Borel subgroup
and $T$ the maximal torus in $B$. The basic facts that are to be
presented in this section can be found e.g. in \cite{LG} Ch. 3,
\cite{Hum} Ch. 1,5, \cite{Bo} Ch. 4, \cite{BGG}. Recall the Levi
decomposition $B=U \cdot T$, where $U$ is the unipotent radical of
$B$. The Weyl group $W$ is by definition $N(T)/T$ where $N(T)$ is
the normalizer of $T$ in $G$. It is generated by the simple
reflections $s_{\beta_1},...,s_{\beta_m}$ corresponding to the
positive simple roots $\Delta= \{\beta_1,...,\beta_m\}$. For $w$
in $W$, the length of $w$, denoted $l(w)$, is the minimum number
$l$ of simple reflections whose product is $w$.
$\Delta^-=\{x_1,...,x_m\}$ denotes the set of negative simple
roots, and $\Phi$ respectively $\Phi^+$ denotes the set of {\em
all} roots respectively the set of all {\em positive} roots. The
longest element of the Weyl group is denoted by $w_0$, and
$B^-=w_0Bw_0$ is the opposite Borel subgroup with its unipotent
radical $U^-=w_0Uw_0$.

Fix also $P$ a parabolic subgroup of $G$ containing $B$. This is
equivalent to choosing $\Delta_P \subset \Delta$, a subset of the
positive simple roots; the Weyl group $W_P$ is generated by the
simple reflections in $\Delta_P$. It is known that each coset
$wW_P$ has a unique representative in $W$ of minimal length. The
set of all such representatives is denoted by $W^P$. The length of
a coset $wW_P$ in $W/W_P$ is by definition the length of its
minimal length representative; the dual $w^\vee$ of $w \in W^P$ is
the minimal length representative of $w_0wW_P$. The codimension
$c(w)$ of $w \in W^P$ is defined to be $l(w^\vee)= \dim G/P -
l(w)$. The dual $s_{\beta}^\vee$ of the simple reflection
associated to $\beta \in \Delta \smallsetminus \Delta_P$ is
denoted by $s(\beta)$ and the minimal length representative
associated to $w_0$ is denoted by $\overline{w}_0$. Thus the
codimension of $s(\beta)$ is equal to $1$ and the codimension of
$\wo$ is zero.

It is known that the cosets $wP$ in $X=G/P$, for $w$ in $W^P$ are
precisely the fixed points of the $T-$action on $X$. Their orbits
$X(w)^o=U \cdot wP$ and $Y(w)^o=U^-\cdot w^\vee P$ are
respectively the Schubert and the opposite Schubert cells,
isomorphic to the affine space $\mathbb{A}^{l(w)}$. The closures
$X(w)=cl(X(w)^o)$ and $Y(w)=cl(Y(w)^o)$ are called the Schubert
varieties (the ordinary and the opposite ones), and they determine
cohomology classes $\sigma(w) = [X(w)]$ respectively
$\tau(w)=[Y(w)]$ in $H^{2c(w)}(X)$ (in fact, because
$Y(w)=w_0X(w)$, the classes $\sigma(w)$ and $\tau(w)$ are equal).
Since the Schubert cells (resp. the opposite cells) cover $X$ by
disjoint affines, their closures $\{\sigma(w) \}$ resp.
$\{\tau(w)\}$ form a $\z-$basis for the cohomology algebra
$H^*(X)$ of $X$ (for these facts, see e.g. \cite{BGG,FW}). Note
that any divisor class in $H^*(X)$ (i.e. a complex degree $1$
class) can be written as an integral linear combination of classes
$\sigma(s(\beta))$ while any curve class (i.e. of complex degree
$\dim X - 1$) can be written as a linear combination of classes
$\sigma(s_\beta)$. In both cases $\beta$ varies over $\Delta
\smallsetminus \Delta_P$.

There is a cohomology pairing \[\langle\textrm{ },\textrm{
}\rangle:H^i(X) \otimes H^j(X) \longrightarrow H^{i+j-2 \dim
X}(pt) \] defined by \[ x \otimes y \longrightarrow \pi_*(x \cup
y) \] where $\pi_*:H^k(X) \longrightarrow H^{k-2 \dim X}(pt)$ is
the Gysin push-forward in cohomology determined by the structure
map $\pi:X \longrightarrow pt$. More about this map is given in
the Appendices from \cite{Mi1,Mi2} and in my thesis. The bases $\{
\sigma(w) \}$ and $\{ \tau(w) \}$ are dual to each other with
respect to this pairing, in the sense that $\langle
\sigma(u),\tau(v)\rangle = \delta_{u,v^\vee}$. This duality
together with its equivariant generalization from \S \ref{eqcohs}
will play an essential role in the proof of the equivariant
quantum Chevalley rule.

A degree $d$ is an integral nonnegative combination $\sum d_\beta
\sigma(s_\beta)$, where the sum is over simple roots $\beta$ in
$\Delta \smallsetminus \Delta_P$. We will often identify the
degree $d$ with the integer sequence $(d_\beta)$. If
$d^{(1)}=(d_\beta^{(1)})$ and $d^{(2)}=(d_\beta^{(2)})$ are two
degrees, we will write $d^{(1)}\geqslant d^{(2)}$ respectively
$d^{(1)}
> d^{(2)}$ if $d_\beta^{(1)}\geqslant d_\beta^{(2)}$ for all $\beta$ respectively
if $d^{(1)}\geqslant d^{(2)}$ but $d^{(1)}$ is not equal to
$d^{(2)}$. To shorten notation, the degree with all components
equal to zero is denoted by $0$. Given $\alpha$ a positive root in
$\pnop$, one defines a degree
\begin{equation}\label{dalpha} d(\alpha) = \sum_{\beta \in
\Delta\smallsetminus \Delta_P} h_\alpha(\omega_\beta)
\sigma(s_\beta) \end{equation} where $h_\alpha$ is the coroot
$2\alpha/(\alpha,\alpha)$; $(\textrm{ } ,\textrm{ } )$ is the
usual inner product on the real subspace of the Lie algebra
$\frak{t}^*$ (the dual of the Lie algebra $\frak{t}$ of $T$); this
subspace is spanned by $\Delta$, and $\omega_\delta$ are the
fundamental weights dual to the basis of simple coroots
$h_\delta$, for $\delta$ in $\Delta$, with respect to the given
inner product. Thus $h_\beta(\omega_\eta)=\delta_{\beta,\eta}$ for
every pair of positive simple roots $\beta$ and $\eta$. For a
geometric interpretation of this degree see \cite{FW} \S 3.

\section{Quantum cohomology of $G/P$}\label{qcohs}

Quantum cohomology of $X=G/P$, denoted $QH^*(X)$, is a graded
$\z[q]-$algebra, where $q$ stands for an indeterminate sequence
$(q_\beta)$, for $\beta$ in $\Delta \smallsetminus \Delta_P$. The
(complex) degree of each $q_\beta$ is $n_{\beta}$, where \[
n_{\beta} = \int_{X(s_\beta)} c_1(TX) \] For a general degree
$d=\sum d_\beta \sigma(s_\beta)$, $q^d$ stands for the monomial
$\prod q_\beta^{d_\beta}$ and $\deg q^d$ denotes $\sum n_\beta
d_\beta$. As a convention, all the degrees will be complex. If
$\alpha$ is a positive root not in $\Phi_P^+$, we denote by
$n(\alpha)$ the degree of $q^{d(\alpha)}$. An explicit formula for
$n(\alpha)$ as well as a geometrical interpretation can be found
in \cite{FW} \S 3. We only note that $n(\alpha)\geqslant 2$ for
any $\alpha$ as above.

The quantum cohomology algebra has a $\z[q]-$basis consisting of
Schubert classes $\sigma(w)$ for $w$ a representative in $W^P$.
The multiplication of two classes $\sigma(u)$ and $\sigma(v)$ is
given by
\[ \sigma(u) \star \sigma(v) = \sum_d \sum_w c_{u,v}^{w,d} q^d
\sigma(w) \] where the first sum is over all sequences
$d=(d_\beta)$ of nonnegative integers (same number of components
as $q$). The coefficients $c_{u,v}^{w,d}$ are the ($3$-pointed,
genus $0$) Gromov-Witten invariants, equal to the number of
rational curves of multidegree $d$ passing through general
translates of Schubert varieties $X(u),X(v)$ and $Y(w^\vee)$.

We need a more formal definition of these coefficients, using
Kontsevich's moduli space of stable maps $\mbotreix$, where
$d=(d_\beta)$ is a multidegree. This is a projective normal
variety of pure (complex) dimension $\dim X + \sum n_{\beta}
d_\beta$ whose closed points consist of stable maps
$f:(C,p_1,p_2,p_3)\longrightarrow X$ of degree $d$ with three
marked points, where $C$ is a tree of $\mathbb{P}^1$'s (see
\cite{FP}, Thm. 2). There are evaluation maps $ev_i:\mbotreix
\longrightarrow X$ sending $(C,p_1,p_2,p_3;f)$ to $f(p_i)$. Then
the coefficient $c_{u,v}^{w,d}$ is given by \[ c_{u,v}^{w,d}=
\pi_*\bigl( ev_1^*(\sigma(u))\cdot ev_2^*(\sigma(v)) \cdot
ev_3^*(\tau(w^\vee)) \bigr) \] in $H^0(pt)$, where $\pi_*$ is the
Gysin push-forward of the structure map $\pi:\mbotreix \rightarrow
pt$ (see \cite{FP,FW} for more details). We also recall the
quantum Chevalley rule proved in \cite{FW}, Thm. 10.1.

\begin{prop}[quantum Chevalley]\label{qc}
Let $\beta$ be a simple root in $\dnop$ and $w$ a minimal length
representative in $W^P$. Then \begin{equation}\label{qcf}
\sigma(s(\beta)) \star \sigma(w) = \sum h_{\alpha}(\omega_\beta)
\sigma(ws_\alpha) + \sum q^{d(\alpha)} h_{\alpha}(\omega_\beta)
\sigma(ws_\alpha)
\end{equation} where the first sum is over all $\alpha \in \pnop$
such that $ws_\alpha$ is a representative in $W^P$ of codimension
$c(ws_\alpha)=c(w)+1$, and the second sum is over those $\alpha
\in \pnop$ such that the $ws_\alpha$ is a representative of
codimension $c(ws_\alpha)=c(w)+1-n(\alpha)$.
\end{prop}

\section{Equivariant cohomology}\label{eqcohs}

The equivariant cohomology of $X$ is the ordinary cohomology of a
``mixed space'' $X_T$, whose definition (see e.g.
\cite{GKM,Br1,Gr,EG} and references therein) we recall. Let $ET
\longrightarrow BT$ be the universal $T-$bundle. The $T-$action on
$X$ induces an action on the product $ET \times X$ by $t\cdot
(e,x)=(et^{-1},tx)$. The quotient space $X_T=(ET \times X)/T$ is
the ``homotopic quotient'' of $X$ and the ($T-$)equivariant
cohomology of $X$ is by definition
\[ H^{i}_T(X)=H^i(X_T). \] In particular, the equivariant
cohomology of a point, denoted by $\Lambda$, is equal to the
ordinary cohomology of the classifying space $BT$. If $\chi$ is a
character in $\widehat{T}=Hom(T, \cx^*)$ it determines a line
bundle \[ L_\chi: ET\times_T \cx_\chi \longrightarrow BT \] where
$\cx_\chi$ is the $1-$dimensional $T-$module determined by $\chi$.
It turns out that the morphism $\widehat{T} \longrightarrow
H^2_T(pt)=\Lambda$ taking the character $\chi$ to the first Chern
class $c_1(L_\chi)$ extends to an isomorphism from the symmetric
algebra of $\widehat{T}$ to $H^*_T(pt)$ (see e.g. \cite{Br1} or \S
7 in \cite{Mi2}). Since the character group is a
finitely-generated free abelian group with basis the negative
simple roots $x_1,...,x_m$, it follows that $\Lambda$ is the
polynomial algebra $\z[x_1,...,x_m]$ in these variables.

The Schubert varieties $X(w)$ respectively $Y(w)$ are stable under
the $T-$action, and determine equivariant cohomology classes
$\sigma(w)^T$ respectively $\tau(w)^T$ in $H^{2c(w)}_T(X)$. Since
$H^*_T(X)$ has a structure of $\Lambda-$algebra (obtained via the
$X-$bundle projection $X_T \longrightarrow BT$), and since the
restriction of the equivariant classes in question to the fiber of
this $X-$bundle determine a basis of the cohomology of the fiber,
the Leray-Hirsch theorem (see \cite{Hus}, Ch. 16) implies that $\{
\sigma(w)^T \}$ respectively $\{ \tau(w)^T \}$ are $\Lambda-$
bases for $H^*_T(X)$.\begin{footnote} {Note that in the
equivariant setting the Schubert classes $\sigma(w)^T$ and
$\tau(w)^T$ are no longer equal; in fact, there is an isomorphism
$\overline{\psi}:H^*_T(X) \to H^*_T(X)$ sending $[X(w)]_T$ to
$[Y(w)]_T$, induced by the involution $\psi:X \to X$ given by
$\psi(x)=w_0\cdot x$. This map is not $T-$equivariant, but it is
equivariant with respect to the map $T \to T$ defined by $t \to
w_0tw_0^{-1}=w_0tw_0$, hence over $H^*_T(pt)$ the isomorphism
$\overline{\psi}$ sends $c_1(L_\chi)$ to
$c_1(L_{w_0\chi})$.}\end{footnote}

These two bases are dual with respect to an equivariant version of
the pairing defined in section \ref{schubert}. This pairing takes
$x \otimes_\Lambda y$ from $H^i_T(X) \otimes_\Lambda H^j_T(X)$ to
$\langle x,y \rangle_T \in H^{i+j-2 \dim X}_T(pt)$ defined by
\[ \langle x,y \rangle_T = \pi_{*}^T(x \cup y) \]
where $\pi_{*}^T:H^k_T(X) \longrightarrow H^{k - 2 \dim X}_T(pt) $
is the equivariant Gysin push-forward (cf. \cite{Mi1,Mi2} or my
thesis). Then the precise formulation of the duality of the two
bases is that
\begin{equation}\label{eqdual} \langle \sigma(u)^T,\tau(v)^T \rangle_T =
\delta_{u,v^\vee}. \end{equation} An equivalent statement is
proved in \cite{Gr}, Lemma 4.2.

The multiplication between two equivariant cohomology classes
$\sigma(u)^T$ and $\tau(v)^T$ is determined by the equivariant LR
coefficients $c_{u,v}^w$:
\[ \sigma(u)^T \cdot \sigma(v)^T = \sum_w c_{u,v}^w \sigma(w)^T \]
These coefficients are homogeneous polynomials in $\Lambda$ of
degree $c(u)+c(v)-c(w)$, and they are equal to the ordinary LR
coefficients if $c(u)+c(v)=c(w)$. Formally, they are defined as \[
c_{u,v}^w = \pi_*^T(\sigma(u)^T \sigma(v)^T \tau(w^\vee)^T) \] and
they turn out to be polynomials in the negative simple roots
$x_1,...,x_m$ with {\em nonnegative} coefficients (this was
conjectured by D. Peterson and has been proved by Graham
\cite{Gr}). A positive formula for them, in this sense, has been
obtained by Knutson and Tao (\cite{KT1}) in the case $X$ is a
Grassmannian. Their proof uses a recurrence formula which is
generalized in this paper (see Cor. \ref{rec} below).

We recall next the formula for the special multiplication with a
divisor class, which will be generalized later. Let $w$ be a
minimal length representative in $W^P$, and let
$w=s_{\beta_{i_1}}\cdot ...\cdot s_{\beta_{i_k}}$ be a reduced
word decomposition. For $\beta_i$ a simple root in $\dnop$, define
the linear form $D(s_{\beta_i},w)$ to be \begin{equation}\label{d}
D(s_{\beta_i},w) = \sum_{i_j=i} s_{\beta_{i_1}}\cdot...\cdot
s_{\beta_{i_{j-1}}}(\beta_i) \end{equation} It can be shown that
each term of the sum is a positive root $\alpha$ in $W$ with the
property that $w^{-1}\alpha$ is a negative root (see the Appendix
for the positivity statement and e.g. \cite{Hum} \S 1.7 for the
second one).

Let $\varphi:\Lambda \longrightarrow \Lambda$ be the automorphism
sending the positive simple root $\beta$ to the negative simple
root $w_0(\beta)$ (see the paragraph before Cor. \ref{properties}
for a proof of that).

\begin{prop}[Equivariant Chevalley formula - see e.g.
\cite{Ku} Thm. 11.1.7(c) and Prop. 11.1.11] \label{eqc} The
following formula holds in $H^*_T(X)$:
\begin{equation}\label{eqcf} \sigma(s(\beta))^T \cdot \sigma(w)^T
= \sum h_{\alpha}(\omega_\beta) \sigma(ws_\alpha)^T +
c_{s(\beta),w}^w \sigma(w)^T  \end{equation} where the sum is over
all $\alpha \in \pnop$ such that $ws_\alpha$ is a representative
in $W^P$ of codimension $c(ws_\alpha)=c(w)+1$, and
$c_{s(\beta),w}^w$ is equal to $\varphi(D(s_\beta,w^\vee))$.
\end{prop}

\section{Equivariant quantum cohomology}\label{eqqcohs}
The equivariant quantum cohomology of $X$, denoted $\eqqhx$, is a
graded $\Lambda[q]-$algebra, where the degree of $q=(q_\beta)$ is
given in \S \ref{qcohs}. It has a $\Lambda[q]-$basis consisting of
Schubert classes $\sigma(w)$, where $w$ varies in $W^P$. The
multiplication, denoted $\circ$, is determined by the 3-pointed
equivariant Gromov-Witten invariants $c_{u,v}^{w,d}$ introduced by
Givental and Kim in \cite{GK}: \[ \sigma(u) \circ \sigma(v) =
\sum_{d \geqslant 0} \sum_w q^d c_{u,v}^{w,d} \sigma(w). \] Recall
that these coefficients are referred to as the equivariant quantum
Littlewood-Richardson coefficients, abbreviated EQLR. They are
homogeneous polynomials in $\Lambda$ of degree $c(u)+c(v)-c(w)-
\sum n_\beta d_\beta$, for $d=(d_\beta)$ (where $n_\beta=\deg
q_\beta$). If the degree $d$ is equal to zero then $c_{u,v}^{w,0}$
is the corresponding equivariant LR coefficient, and if the
polynomial degree is equal to zero then $c_{u,v}^{w,d}$ is the
quantum LR coefficient. The coefficients for which both $d>0$ and
the polynomial degree is larger then zero will be called {\bf
mixed}.

The formal definition of the EQLR coefficients is similar to the
one of the quantum ones, except that all the maps and cohomology
classes are replaced by their equivariant versions. More
precisely, the $T-$action on $X$ induces an action on the moduli
space of stable maps $\mbotreix$ given by:
\[ t \cdot (C,p_1,p_2,p_3;f)
:= (C,p_1,p_2,p_3;\tilde{f}) \] where $\tilde{f}(x):=t \cdot
f(x)$, for $x$ in $C$ and $t$ in $T$. The evaluation maps
$ev_i:\mbotreix \longrightarrow X$ ($i=1,2,3$) and the structure
map to the point $\pi:\mbotreix \longrightarrow pt$ are
$T-$equivariant. Following \cite{Kim2} \S 3.1 define the
equivariant Gromov-Witten invariant
\begin{equation}\label{eqlrdef} c_{u,v}^{w,d} =
\pi_*^T\bigl((ev_1^T)^* (\sigma(u)^T) \cdot (ev_2^T)^*
(\sigma(v)^T) \cdot (ev_3^T)^* (\tau(w^\vee)^T)\bigr)
\end{equation} where $\pi_*^T:H^k_T(\mbotreix) \longrightarrow
H^{k - 2 \dim \mbotreix}_T(pt)$ is the equivariant Gysin push
forward. More details about this definition, as well as proofs of
the properties of the EQLR coefficients can be found in
\cite{Mi1,Mi2}.

\section{Equivariant quantum Chevalley rule}\label{eqcs}

The aim of this section is to prove the equivariant quantum
Chevalley rule. As in \cite{Mi1}, this will follow from a certain
vanishing property of the EQLR coefficients. We need the following
lemma:
\begin{lemma}\label{vanishing} Let $u,v,w$ be three
representatives in $W^P$ such that one of the intersections
$ev_1^{-1}(X(u)) \cap ev_3^{-1}(Y(w^\vee))$ or $ev_2^{-1}(X(v))
\cap ev_3^{-1}(Y(w^\vee))$ in $\mbotreix$ is empty. Then the EQLR
coefficient $c_{u,v}^{w,d}$ is equal to zero.
\end{lemma}
\begin{proof} The hypothesis implies that one of the products
$ev^{T*}_1(\sigma(u)^T) \cdot ev^{T *}_3(\tau(w^\vee)^T)$ or
$ev^{T *}_2(\sigma(v)^T) \cdot ev^{T *}_3(\tau(w^\vee)^T)$ in
$H^*_T(X)$ vanishes (for details see the Fact 1 in the proof of
Lemma 4.3, \cite{Mi1}). The assertion follows then from the
definition of the EQLR coefficients (see formula (\ref{eqlrdef})
above).
\end{proof}

The following Lemma, inspired from \cite{FW}, uses a weaker
version of Kleiman's transversality Theorem (\cite{Kl}) to show
that the intersection of the inverse images of two opposite
Schubert varieties through a $G-$equivariant map has the expected
dimension. This is in fact the key lemma used in the proof of the
vanishing of the mixed EQLR coefficients from the EQ Chevalley
formula (see Lemma \ref{main} below).

\begin{lemma}\label{Kleiman}Let $Z$ be a reduced,
possibly reducible, pure dimensional $G$-variety and let $F:Z
\longrightarrow X \times X$ be a $G-$equivariant morphism, where
$G$ acts diagonally on $X \times X$. Then, for any $u$ and $v$ in
$W^P$, the subscheme $F^{-1}(X(u) \times Y(v))$ is either empty or
of codimension $c(u)+c(v)$.
\end{lemma}
\begin{proof} Note that every irreducible component of $Z$ must
be $G-$invariant. Indeed, the action of $G$ permutes the
irreducible components, and the unit element in $G$ fixes each
component. This implies that it is enough to show the result in
the case when $Z$ is irreducible, which is a part of what is
proved in Lemma 7.2 from \cite{FW}. For convenience, we summarize
its proof. Kleiman's transversality result (\cite{Kl}, Thm. 2)
yields an open subset $U$ in $G \times G$, invariant under the
diagonal left-multiplication by $G$, such that $F^{-1}(h_1X(u)
\times h_2X(v))$ is either empty or of codimension $c(u)+c(v)$ for
any $(h_1,h_2)$ in $U$. Another lemma (\cite{FW}, Lemma 7.1),
shows that for any pair $(g_1,g_2)$ in $G \times G$ such that the
intersection $g_1Bg_1^{-1} \cap g_2Bg_2^{-1}$ is a maximal torus
in $G$, there is a pair $(h_1,h_2)$ in $U$ such that
\begin{center} $h_1 X(u) = g_1X(u)$ and  $h_2 X(v) = g_2X(v)$.
\end{center} The result follows then by taking $g_1=1$ and
$g_2=w_0$.
\end{proof}

{\it Remark:} Since only the dimension assertion of the Kleiman's
Theorem is used, the lemma is valid in all characteristics.
However, it will be used for $Z$ being the moduli space of stable
maps $\mbotreix$, whose construction is done in characteristic
zero.\\

We are ready to prove the main vanishing result for the EQLR
coefficients. Recall that $n_{\beta}$ denotes the complex degree
of the indeterminate $q_\beta$, for $\beta$ in $\dnop$.

\begin{lemma}[Main Lemma]\label{main} Let $u,v,w$ be
representatives in $W^P$ and $d=(d_\beta)$ a nonzero degree such
that $c(u)+1 > c(w)+ \sum n_{\beta} d_\beta$. Then the EQLR
coefficient $c_{u,v}^{w,d}$ is equal to zero.
\end{lemma}
\begin{proof} By Lemma \ref{vanishing} it is enough to
show that the intersection $E(u,v):=ev_1^{-1}(X(u)) \cap
ev_3^{-1}(Y(w^\vee))$ in $\mbotreix$ is empty. The hypothesis
implies that \[ c(u) + c(w^\vee)+1 > \dim X + \sum
n_{\beta}d_\beta = \dim \mbotreix
\] which is equivalent to \[ c(u)+ c(w^\vee) \geqslant \dim
\mbotreix . \] Lemma \ref{Kleiman} applied to $Z=\mbotreix$ and
$F:\mbotreix \longrightarrow X \times X$ given by $F=(ev_1, ev_3)$
implies that the intersection in question is at most finite.
Moreover, the boundary $B$ of $\mbotreix$, which is the subvariety
consisting of stable maps $(C,p_1,p_2,p_3;f)$ where the curve $C$
has at least one node is $G-$invariant and it is of codimension
one. Applying again Lemma \ref{Kleiman} for $Z=B$ and $F$
restricted to $B$, shows that if $E(u,v)$ is not empty then it
cannot intersect $B$, so all its points must be stable maps whose
sources are curves isomorphic to $\mathbb{P}^1$.

One the other side, given a stable map $f:(C,p_1,p_2,p_3)
\longrightarrow X$ in $\mbotreix$ such that $f(p_1)$ is in $X(u)$
and $f(p_3)$ is in $Y(w^\vee)$ and $C \simeq \mathbb{P}^1$, one
can produce a curve in $E(u,v)$ by letting $f(p_2)$ to vary in the
image of $C$ through $f$ (this image is not a point, since the
degree $d$ is not zero). This constitutes a contradiction with the
fact that $E(u,v)$ is finite.
\end{proof}

{\em Remarks:} 1. When $X$ is a Grassmannian, this vanishing
result is weaker than the one obtained in \cite{Mi1}.

2. It is known that the moduli space $\mbotreix$ is irreducible
(\cite{KP,T}), so we could have used to original version of Lemma
\ref{Kleiman}, where $Z$ was assumed irreducible.\\

An immediate consequence of the Main Lemma is that all the mixed
EQLR coefficients in the product $\sigma(s(\beta))\circ \sigma(w)$
must vanish. Indeed, by definition, a coefficient
$c_{s(\beta),u}^{w,d}$ is mixed if the degree $d=(d_\beta)$ is not
zero and if its polynomial degree is also not zero, i.e. if
$c(u)+1 > c(w) + \sum n_{\beta}d_\beta$. This is precisely the
hypothesis of the Main Lemma, so $c_{s(\beta),u}^{w,d}=0$. In
particular, combining the quantum and equivariant Chevalley
formulae (\ref{qcf}) and (\ref{eqcf}) yields

\begin{thm}[Equivariant quantum Chevalley rule]\label{eqqc} Let
$\beta$ be a simple root in $\dnop$ and $w$ a minimal length
representative in $W^P$. Then the following formula holds in the
equivariant quantum cohomology of $X$
\begin{equation}\label{eqqcf}
\sigma(s(\beta)) \circ \sigma(w) = \sum h_{\alpha}(\omega_\beta)
\sigma(ws_\alpha) + \sum q^{d(\alpha)} h_{\alpha}(\omega_\beta)
\sigma(ws_\alpha) + c_{s(\beta),w}^w \sigma(w)
\end{equation} where the first sum is over all positive roots
$\alpha$ in $\pnop$ such that $ws_\alpha$ is a representative in
$W^P$ and $c(ws_\alpha)=c(w)+1$, and the second sum is over those
$\alpha$ in $\pnop$ such that $ws_\alpha$ is a representative in
$W^P$ of codimension $c(ws_\alpha)=c(w)+1 - n(\alpha)$. The
coefficient $c_{s(\beta),w}^w$ is the one given in Prop.
\ref{eqc}. \end{thm}

As in the Grassmannian case, the equivariant quantum Chevalley
formula implies a recursive formula satisfied by the EQLR
coefficients. Using double recursion, first on the degree $d$,
then on the polynomial degree, this formula shows that the EQLR
coefficient $c_{u,v}^{w,d}$ is equal to a homogeneous combination
of EQLR coefficients, some with smaller degree $d$, and the
remaining ones with the same degree $d$ but {\em higher}
polynomial degree.

\begin{cor}\label{rec} Let $u,v,w$ be representatives in $W^P$,
$\beta$ a positive simple root in $\dnop$ and $d$ a degree. Then
the following formula holds:
\begin{eqnarray*}\label{recf}
(c_{s(\beta),w}^w - c_{s(\beta),u}^u) c_{u,v}^{w,d} & = &
\sum_\alpha h_\alpha(\omega_\beta) c_{u_1,v}^{w,d} - \sum_\alpha
h_\alpha(\omega_\beta) c_{u,v}^{w_1,d} + \\ & & \sum_{\alpha, d
\geqslant d(\alpha)} h_\alpha(\omega_\beta)
c_{u_2,v}^{w,d-d(\alpha)} - \sum_{\alpha,d\geqslant d(\alpha)}
h_\alpha(\omega_\beta) c_{u,v}^{w_2,d-d(\alpha)} \end{eqnarray*}
where the first two sums are over $\alpha \in \pnop$ such that
$u_1=us_\alpha$ is a representative in $W^P$ of codimension
$c(u_1)=c(u)+1$ and $w_1$ is such that $w_1s_\alpha=w$ and
$c(w)=c(w_1)+1$; the last two sums are over $\alpha \in \pnop$
with $d \geqslant d(\alpha)$ (with $d(\alpha)$ defined by
(\ref{dalpha})) such that $u_2=us_\alpha$ is a representative in
$W^P$ with $c(u_1)=c(u)+1-n(\alpha)$ and $w_2$ is a representative
such that $w=w_2s_\alpha$ and $c(w)=c(w_2)+1-n(\alpha)$.\end{cor}
\begin{proof} This is a straightforward computation. As in the Grassmannian
case (see \cite{Mi1}, Prop. 5.1), the formula follows by
collecting the coefficient of $q^d \sigma(w)$ from both sides of
the associativity equation
\[ \bigl(\sigma(s(\beta))\circ \sigma(u)\bigr)\circ \sigma(v) =
\sigma(s(\beta))\circ \bigl( \sigma(u)\circ \sigma(v)\bigr).
\] \end{proof}

{\em Remark:} This formula is the main ingredient for an effective
algorithm to compute the EQLR coefficients, found in \S  \ref{algs} below.

\section{Two formulae}\label{formulaes}

The aim of this section is to prove two formulae satisfied by the
EQLR coefficients, which will be used in the algorithm of the next
section.  From now on, all the results will be algorithmic, and a
coefficient $c_{u,v}^{w,d}$ is regarded as a (possibly rational)
homogeneous function of degree $c(u)+c(v)-c(w)-\sum
n_{\beta}d_\beta$. The latter quantity will still be called {\em
polynomial degree}, even though it may a priori be negative.

The formulae to be proved are the natural generalizations to
$X=G/P$ of the formulae from Propositions 6.1 and 6.2 in
\cite{Mi1} (where $X$ was a Grassmannian). To state them, we
define a reversed Bruhat ordering denoted $\prec$ on the
permutations $W^P$ as follows: write $w_1 \rightarrow w_2$ if
there exists $\alpha$ a positive root in $\pnop$ such that
$w_2=w_1s_\alpha$ and $c(w_2)>c(w_1)$. Then $u \prec w$ if there
is a chain $u=w_0 \rightarrow w_1 \rightarrow ... \rightarrow
w_k=w$ (in fact, in the definition of $w_1\rightarrow w_2$, \S
5.11 in \cite{Hum} shows that it is enough to consider those
$\alpha$ for which $c(w_2)=c(w_1)+1$).

We rewrite next the recursive formula from Cor. \ref{rec} as
\begin{equation}\label{mrecf} \begin{split} c_{u,v}^{w,d} & =
\frac{\sum_\alpha h_\alpha(\omega_\beta)
c_{u_1,v}^{w,d}}{F_{w,u}(\beta)}
- \frac{\sum_\alpha h_\alpha(\omega_\beta) c_{u,v}^{w_1,d}}{F_{w,u}(\beta)} + \\
&  \frac{\sum_{\alpha, d \geqslant d(\alpha)}
h_\alpha(\omega_\beta) c_{u_2,v}^{w,d-d(\alpha)}}{F_{w,u}(\beta)}
- \frac{\sum_{\alpha,d\geqslant d(\alpha)} h_\alpha(\omega_\beta)
c_{u,v}^{w_2,d-d(\alpha)}}{F_{w,u}(\beta)} \end{split}
\end{equation} where $w$ is different from $u$ and $F_{w,u}(\beta)$
is the linear homogeneous
form in the negative simple roots in $W$ defined by
\[ F_{w,u}(\beta)= c_{s(\beta),w}^w-c_{s(\beta),u}^u \]
with $\beta$ a positive simple root in $\dnop$.

The following Lemma gives some of the basic properties of the
forms $F_{w,u}(\beta)$ together with the equivariant coefficients
$c_{s(\beta),w}^w$ used to define it. To shorten notations, for
each $u,w$ in $W^P$ such that $u \prec w$, we denote by $Cov(u,w)$
the set of positive roots $\alpha$ in $\pnop$ such that
$us_\alpha$ is a representative in $W^P$, $c(us_\alpha)=c(u)+1$
and $us_\alpha \prec w$.
\begin{lemma}\label{F} Let $u,w$ be two distinct representatives in $W^P$.
\begin{enumerate} \item The coefficient $c_{s(\beta),w}^w$, for $\beta$ in $\dnop$,
is a linear homogeneous combination of negative simple roots
$x_1,...,x_m$ with nonnegative coefficients, and there exists a
$\beta$ for which this coefficient is nonzero.

\item There exists a positive simple root $\beta$ in $\dnop$ such
that $F_{w,u}(\beta)$ is nonzero.

\item Assume that $u \prec w$ in the reverse Bruhat ordering
previously defined. Then for any $\beta$ in $\dnop$ the form
$F_{w,u}(\beta)$ is a linear nonnegative combination of negative
simple roots.

\item If $u \prec w$ the set $Cov(u,w)$ is nonempty. Moreover, for
any positive root $\alpha$ in $Cov(u,w)$ and any $\beta$ positive
simple root in $\dnop$ such that $F_{us_\alpha,u}(\beta)$ is
nonzero, the integer $h_\alpha(\omega_\beta)$ is positive.
\end{enumerate}
\end{lemma}
\begin{proof} This is Cor. \ref{properties} from the Appendix.
\end{proof}

The next proposition shows that the computation of the EQLR
coefficients $c_{u,v}^{w,d}$ such that $u$ is not less than $w$
with respect to the ordering $\prec$ can be reduced to the
computation of some coefficients with smaller degree $d$.

\begin{prop} Let $u,v,w$ be three representatives in $W^P$ and $d=(d_\beta)$
a degree. Assume that $u \nprec w$. Then
\begin{equation}\label{inclusion} c_{u,v}^{w,d}=E_{u,v,w}(d)
\end{equation} where $E_{u,v,w}(d)$ is algorithmically known, and it is
a linear homogeneous form in EQLR coefficients of degree strictly
smaller than $d$, with coefficients in the fraction field
$R(\Lambda)$ of $\Lambda=\z[x_1,...,x_m]$. If $d=0$ then
$E_{u,v,w}(0)=0$.
\end{prop}
\begin{proof} We argue by descending induction on
$c(u)-c(w) \leqslant \dim X$, with equality exactly when $u$ is
the unit element $1$ in $W^P$ and $w=\overline{w}_0$ (the
representative in $W^P$ which indexes the unit element in
$\eqqhx$). In this case the first two sums from (\ref{mrecf})
(applied with a suitable $\beta$, such that
$F_{\overline{w}_0,1}(\beta)$ doesn't vanish, cf. Lemma \ref{F})
disappear and hence $c_{1,v}^{\overline{w}_0,d}$ is equal to a
combination of EQLR coefficients coming from the last two sums,
which have coefficients of degree strictly smaller than $d$ (thus
defining $E_{1,v,\overline{w}_0}(d)$). Assume now that
$c(u)-c(w)<\dim X$. Since $u \nprec w$, $u$ cannot be equal to
$w$, hence one can apply formula (\ref{mrecf}) to $c_{u,v}^{w,d}$
(again, with a suitable $\beta$). The last two sums enter into the
definition of $E_{u,v,w}(d)$. Note that the coefficients
$c_{u_1,v}^{w_1,d}$ from the first two sums satisfy $c(u)-c(v) <
c(u_1)-c(w_1)$, so to finish the proof it is enough to show that
these coefficients satisfy the induction hypothesis, i.e. that
$u_1 \nprec w_1$. Assuming the contrary, i.e. $u_1 \prec w_1$,
since $u \prec u_1$ and $w_1 \prec w$, it follows that $u \prec
w$, a contradiction. The case $d=0$ is treated in the same way,
proving now that $c_{u,v}^{w,0}=0$ if $u \nprec w$.
\end{proof}
Next is a formula which shows that the coefficients of the form
$c_{w,w}^{w,d}$ are determined algorithmically by coefficients of
the form $c_{u_0,w}^{w,d}$ with $u_0 \prec w$ and by coefficients
of (strictly) smaller degree $d$. To write this formula we will
introduce a weighted Bruhat-type oriented graph $\Gamma(u,w;f)$,
encoding all the possible saturated paths $\pi$ from $u$ to $w
\succ u$, weighted according to the recipe described below.

\begin{defn}\label{dg} Let $u,w$ be two representatives in $W^P$ such that $u \prec
w$. The weighted oriented graph $\Gamma(u,w;f)$ is given by the
following data:

\begin{itemize} \item A set $V(u,w)$, of vertices, which consists of all
representatives $v$ in $W^P$ such that $u \prec v \prec w$.

\item For each $v_1,v_2$ in $V(u,w)$ such that $v_2=v_1s_\alpha$,
where $\alpha$ is in $\pnop$ such that $c(v_2)=c(v_1)+1$ (i.e.
$\alpha$ is in $Cov(v_1,w)$), there is an edge between $v_1$ and
$v_2$, oriented from $v_1$ to $v_2$. This edge is denoted in short
by $v_1 \xrightarrow{\alpha} v_2$.

\item An assignment $f:V(u,w) \longrightarrow \dnop$, $v
\rightarrow \beta(v)$ such that $F_{w,v}(\beta(v))$ is a nonzero,
nonnegative combination of negative simple roots (such an $f$
exists by Lemma \ref{F}, assertions (2) and (3)).

\end{itemize}

For each edge $v \xrightarrow{\alpha}vs_\alpha$ in $\Gamma(u,w;f)$
define its weight to be

\[wt(v \xrightarrow{\alpha}vs_\alpha)=
\frac{h_\alpha(\omega_{\beta(v)})}{F_{w,v}(\beta(v))}
\]

A {\bf path} $\pi$ in $\g$ is any oriented path from $u$ to $w$.
The weight of $\pi$, denoted $wt(\pi)$, is the product of all the
weights of the edges it contains. \end{defn}

To any such oriented weighted graph $\g$ one associates a
homogeneous rational function in $R(w,u;f)$ in the fraction field
$R(\Lambda)$ defined by

\begin{equation}\label{r} R(u,w;f) = \sum_\pi wt(\pi)
\end{equation}
where the sum is over all paths in the graph $\g$. We will need
the following Lemma:

\begin{lemma}\label{noz} Let $u,w$ two representatives in $W^P$ such that
$u \prec w$. Then there exists an assignment $f:V(u,w)
\longrightarrow \dnop$ as in Definition \ref{dg} above such that
the function $R(u,w;f)$ is not zero. \end{lemma}
\begin{proof} For any assignment $f$, a path
\[ \pi: u=u_0 \xrightarrow{\alpha_1} u_1
\xrightarrow{}...\xrightarrow{\alpha_k}u_k=w \] in $\g$ has weight
\[
wt(\pi)=\prod_{i=1}^{k}\frac{h_{\alpha_{i}}\bigr(\omega_{\beta(u_{i-1})}\bigl)
}{F_{w,u_{i-1}}(\beta(u_{i-1}))}\] Definition \ref{dg} and Lemma
\ref{F} imply that the denominator is a nonzero, nonnegative
homogeneous linear combinations of negative simple roots, while
the numerator is a product of nonnegative integers. Then an
assignment for which $R(u,w;f)$ is nonzero would be any assignment
for which the graph $\g$ has a path $\pi$ of nonzero weight. Such
an assignment is provided by (4) from Lemma \ref{F}. To each
vertex $v$ from $V(u,w)$, choose a positive root $\alpha(v)$ in
$Cov(v,w)$. Then $\beta(v)$ is chosen to be any of the positive
simple roots in $\dnop$ such that $F_{vs_{\alpha(v)},v}(\beta(v))$
is nonzero. Indeed, in this case
$F_{w,v}(\beta(v))=F_{w,vs_{\alpha(v)}}(\beta(v))+
F_{vs_{\alpha(v)},v}(\beta(v))>0$ by Lemma \ref{F} and by
definition of $\beta(v)$. Given this assignment, a path $\pi$ of
nonzero weight can be constructed inductively as follows: assuming
that the $(i+1)^{st}$ element $u_i$ is constructed, and a positive
root $\alpha=\alpha(u_i)$ has been chosen as before, $u_{i+1}$ is
equal to $u_is_{\alpha}$. In this case, the weight of the edge
$u_i \xrightarrow{\alpha} u_{i+1}$ is
\[\frac{h_{\alpha}(\omega_{\beta(u_i)})}{F_{w,u_i}(\beta(u_i))}\]
and $h_{\alpha}(\omega_{\beta(u_i)})$ is a positive integer by
Lemma \ref{F}(4), as desired.
\end{proof}

\noindent {\it Remark:} We will see (Cor. \ref{assign} below) that
the function $R(u,w;f)$ is independent of the assignment $f$, and
hence nonzero, by the Lemma just proved. This will be a
consequence of the fact that $c_{\overline{w}_0,w}^{w,0}=1$ and of
the formula (\ref{equalityf}), which is proved next.

\begin{prop}\label{equality} Assume that the EQLR coefficients are
commutative, i.e. $c_{u,v}^{w,d}=c_{v,u}^{w,d}$ for any
representatives $u,v,w$ in $W^P$ and let $u_0,w$ be two such
representatives such that $u_0 \prec w$. Then for all assignments
$f:V(u_0,w) \longrightarrow \dnop$, the EQLR coefficient
$c_{u_0,w}^{w,d}$ satisfies the following formula
\begin{equation}\label{equalityf} c_{u_0,w}^{w,d} =
R(u_0,w;f)c_{w,w}^{w,d} + E'(u_0,w;f,d)
\end{equation} where $E'(u_0,w;f,d)$ is algorithmically known, and it is an
$R(\Lambda)-$linear homogeneous expression in coefficients of
degree strictly less than $d$. If $d=0$ then $E'(u_0,w;f,0)=0$.
Moreover, there exists an assignment $f$ such that $R(u_0,w;f)$ is
nonzero. \end{prop}
\begin{proof} We use ascending induction on $c(w)-c(u_0) \geqslant 0$.
The base of the induction, i.e $w=u_0$, is obvious. In this case
$E'(w,w;f,d)=0$ and $R(w,w;f)=1$. Assume that $c(w)-c(u_0)>0$.
Applying (\ref{mrecf}) to $c_{u_0,w}^{w,d}$ and using Prop.
\ref{inclusion} yields
\begin{equation}\label{indeq} c_{u_0,w}^{w,d}  = \bigl(\sum_\alpha
h_\alpha(\omega_{\beta(u_0)})
c_{u_0s_\alpha,w}^{w,d}\bigr)/F_{w,u_0}(\beta(u_0)) +
E'_{w,u_0}(\beta(u_0),d) \end{equation} where the sum is over all
$\alpha$ in $Cov(u_0,w)$ and $E'_{w,u_0}(\beta(u_0),d)$ is an
$R(\Lambda)$-linear combination of the EQLR coefficients of
strictly smaller degree. It contains the EQLR coefficients from
the last two sums of (\ref{mrecf}) since they have smaller degree,
the coefficients $c_{u_0,w}^{w_1,d}$ from the second sum, since
they satisfy $c_{u_0,w}^{w_1,d}=E_{u_0,w,w_1}(d)$ by Prop.
\ref{inclusion} (because $c_{u_0,w}^{w_1,d}=c_{w,u_0}^{w_1,d}$ and
$w \nprec w_1$), and those coefficients $c_{u_0s_\alpha,w}^{w,d}$
from the first sum for which $\alpha$ is in $\pnop$,
$c(u_0s_\alpha)=c(u_0)+1$ but $u_0s_\alpha \nprec w$, to which one
applies again Prop. \ref{inclusion}. Note that (\ref{indeq}) is
equivalent to
\begin{equation}\label{indeq2}
c_{u_0,w}^{w,d}= \sum_{\alpha \in Cov(u_0,w)}
wt(u_0\xrightarrow{\alpha}u_0s_\alpha) c_{u_0s_\alpha,w}^{w,d} +
E'_{w,u_0}(\beta(u_0),d)
\end{equation} Induction hypothesis, applied to each
$c_{u_0s_\alpha,w}^{w,d}$ (for $\alpha$ in $Cov(u_0,w)$), implies
that
\begin{equation}\label{indeq3} c_{u_0s_\alpha,w}^{w,d} = R(u_0s_\alpha,w;f_\alpha)c_{w,w}^{w,d} +
E'(u_0s_\alpha,w;f_\alpha,d)\end{equation} where $f_\alpha$ is the
restriction $f_{|V(u_0s_\alpha,w)}:V(u_0s_\alpha,w)
\longrightarrow \dnop$. Combining (\ref{indeq2}) and
(\ref{indeq3}) implies that modulo the coefficients of degree
smaller than $d$, $c_{u_0,w}^{w,d}$ is equal to
\begin{equation}\label{indeq4} c_{u_0,w}^{w,d}= \sum_{\alpha \in
Cov(u_0,w)} wt(u_0\xrightarrow{\alpha}u_0s_\alpha)
R(u_0s_\alpha,w;f_\alpha)c_{w,w}^{w,d}\end{equation} But the
coefficient of $c_{w,w}^{w,d}$ on the right is precisely
$R(u_0,w;f)$, which proves the first assertion of the Proposition.

For the second assertion, by Prop. \ref{inclusion}, all the
coefficients of the form $c_{u_1,v_1}^{w_1,0}$ vanish if either
$u_1 \nprec w_1$ or $v_1 \nprec w_1$. Retracing the procedure used
in the proof yields that
\[ c_{u_0,w}^{w,0}=R(u_0,w;f)c_{w,w}^{w,0} \] which proves the
second assertion of the Proposition. The last one is an immediate
consequence of Lemma \ref{noz}. \end{proof}

\section{An algorithm to compute the EQLR coefficients}\label{algs}
As in the Grassmannian case (\cite{Mi1}), the algorithm to compute
the EQLR coefficients is by double induction: on the degree $d$ of
a coefficient $c_{u,v}^{w,d}$, then descending induction on the
polynomial degree $c(u)+c(v)-c(w)- \sum_{\beta \in \dnop} n_\beta
d_\beta$. Its main ingredient is the formula (\ref{mrecf}), which
writes the EQLR coefficient $c_{u,v}^{w,d}$ as a combination of
coefficients with smaller degree $d$ and coefficients of the same
degree $d$, but with larger polynomial degree. Assuming
commutativity of the EQLR coefficients, applying (\ref{mrecf})
repeatedly reduces the computation of any coefficient to
coefficients of smaller degree and coefficients of the form
$c_{\eta,\eta}^{\eta,d}$. The latter ones can be computed using
the formula (\ref{equalityf}), for $u_0=\overline{w}_0$ and
$w=\eta$ (note that $\overline{w}_0$ is the smallest element in
$W^P$ with respect to the ordering $\prec$, and indexes the unit
in $\eqqhx$). Recall the assumption made at the beginning of \S
\ref{formulaes} that the EQLR coefficients are homogeneous
rational functions of the expected degree.
\begin{thm}\label{alg} The EQLR coefficients are determined algorithmically
by the following formulae:
\begin{enumerate} \item[(a)] (multiplication by unit)
$$ c_{\overline{w}_0,w}^{w,d}= \left \{
\begin{array}{ll} 1 & \textrm{if } d=0\\
0 & \textrm { } \textrm{otherwise}
\end{array} \right.$$

\item[(b)] (commutativity) $c_{u,v}^{w,d}=c_{v,u}^{w,d}$ for any
$u,v,w$ in $W^P$ and any degree $d$.

\item[(c)] The formula (\ref{mrecf}), for $u \neq w$.
\end{enumerate} \end{thm}

Before proving the Theorem, we present an immediate consequence of
it:
\begin{cor}\label{coralg} Let $A$ be a graded, commutative,
associative $\Lambda[q]-$algebra with unit, where the degree of
$q=(q_\beta)_{\beta \in \dnop}$ is as usual (see \S \ref{qcohs}).
Assume that:

\noindent 1. $A$ has an additive $\Lambda[q]-$basis $\{t_w\}_{w
\in W^P}$ (graded as usual).

\noindent 2. The equivariant quantum Chevalley formula
(\ref{eqqcf}) holds.

Then $A$ is canonically isomorphic to $\eqqhx$, as
$\Lambda[q]-$algebras.
\end{cor}

\begin{proof}[Proof of the Corollary] The structure constants of
$A$ clearly satisfy (a) and (b); (c) follows from the
associativity of $A$ and the EQ Chevalley formula (cf. Cor.
\ref{rec}).
\end{proof}

\begin{proof}[Proof of Thm. \ref{alg}] The proof is by double
induction: ascending induction on the degree $d$ and descending
induction on the polynomial degree. For each fixed $d$ there are
two main steps in the algorithm: the first is to compute the
coefficients of the form $c_{w,w}^{w,d}$, using formula
(\ref{equalityf}), and the second is to (algorithmically) compute
all the other coefficients. The second step uses descending
induction on the polynomial degree.

The base of the induction on $d$ is the situation when $d=0$.

{\it Step 1: Compute the coefficients of the form
$c_{w,w}^{w,0}$.}

We use formula (\ref{equalityf}) with $u_0=\overline{w}_0$. By
(a), the coefficient $c_{\overline{w}_0,w}^{w,0}$ is equal to $1$,
for any $w$ in $W^P$. Choose an assignment $f:V(\overline{w}_0,w)
\longrightarrow \dnop$ such that $R(\overline{w}_0,w;f)$ is not
equal to zero (this assignment exists by Lemma. \ref{noz}). Then
\[ c_{w,w}^{w,0} = 1/R(\overline{w}_0,w;f). \]

{\it Step 2: Compute all other coefficients of degree $0$.} We
argue by descending induction on the polynomial degree. Note that
the last two sums of (\ref{mrecf}) vanish since $d=0$ and
$d(\alpha)>0$. The largest polynomial degree for the coefficient
$c_{u,v}^{w,0}$ is achieved exactly when $u=v=1$ and
$w=\overline{w}_0$. In this case the coefficient in question is
equal to $0$, by Prop. \ref{inclusion}, since $1 \nprec
\overline{w}_0$. Consider now a coefficient $c_{u,v}^{w,0}$ of
smaller polynomial degree. If $u=v=w$ this coefficient is known by
Step 1. By commutativity, we can assume then that $u \neq w$.
Applying formula (\ref{mrecf}) to $c_{u,v}^{w,0}$ (with a suitable
$\beta$) writes this coefficient as an $R(\Lambda)-$linear
combination of coefficients of polynomial degree larger by one,
hence known by induction. This finishes the proof of the case
$d=0$.

Assume now that $d$ is not zero. Note that the induction on $d$
allows us to ignore all the terms of degree less than $d$ in
equations (\ref{mrecf}),(\ref{inclusion}) and (\ref{equalityf}).
But then the proof is the same as the one to the base case, since
the equations obtained in the previous manner are the same as
those for the case $d=0$. \end{proof}

\subsection{Remarks:} 1. The algorithm of Thm. \ref{alg} provides in
particular an algorithm to compute the 3-pointed Gromov-Witten
coefficients for {\it any} homogeneous variety $G/P$. Different
algorithms for that are discussed in the introduction.

2. The hypothesis (a) in the Theorem can be changed to

(a') The coefficient
$c_{\overline{w}_0,\overline{w}_0}^{\overline{w}_0,d}$ is equal to
$0$ unless $d=0$ when it is equal to $1$.

\noindent  and

(a'') The equivariant quantum Chevalley terms
$c_{s(\beta),v}^{w,d}$, as given in the formula (\ref{eqqcf}).

The proof goes the same way, except that in Step 1 the coefficient
$c_{w,w}^{w,d}$ is computed starting from $c_{s(\beta),w}^{w,d}$
with $s(\beta) \prec w$. Note that this requires less number of
computations, but one has to input the equivariant quantum
Chevalley coefficients, which are known anyway from formula
(\ref{eqqcf}).

3. To reduce further the number of computations needed in the
algorithm, one can also impose the following conditions:
\begin{enumerate}
\item[(i)] $c_{u,v}^{w,d}$ is equal to $0$ if it has negative
polynomial degree.

\item[(ii)] $c_{u,v}^{w,d}$ is equal to $0$ if $c(u)+1 > c(w)+
\sum n_{\beta} d_\beta$ or  $c(v)+1 > c(w)+ \sum n_{\beta}
d_\beta$ by Lemma \ref{main}.
\end{enumerate}

4. The number of computations can be reduced even further than in
previous remark, provided one knows some pure quantum coefficient
of the form $c_{u,w}^{w,d}$ with $u \prec w$. As in Remark 2, this
coefficient is then used in Step 1 to compute $c_{w,w}^{w,d}$.

\section{Consequences in equivariant cohomology of $G/P$.}
\label{eqconseq} A by-product of the previous algorithm (more
precisely of the formulae (\ref{inclusion}) and (\ref{equalityf}))
is a formula for the equivariant coefficients of the form
$c_{u,w}^w=c_{u,w}^{w,0}$. Formulae for these coefficients were
known before (\cite{Bi}) therefore we obtain in particular some
interesting combinatorial identities.

If $u \nprec w$ formula (\ref{inclusion}) implies that
$c_{u,w}^w=0$. For $u \prec w$ formula (\ref{equalityf}) together
with the fact that $c_{\overline{w}_0,w}^w=1$ shows that
\begin{equation}\label{top} c_{w,w}^w = 1/R(\overline{w}_0,w;f)
\end{equation} for any assignment $f$ (cf. Def. \ref{dg}). If $u
\prec w$, but $u$ is not equal to $w$, using the same formula
yields
\begin{equation}\label{uw} c_{u,w}^w=R(u,w;f')/R(\overline{w}_0,w;f)
\end{equation} for any assignments $f$ and $f'$. We recall a
different formula for the coefficient $c_{u,w}^w$, proved in
\cite{Bi} (see also Prop. 11.1.11 from \cite{Ku})\begin{footnote}{
The equivariant coefficient $c_{u,w}^w$ is equal to
$\varphi\bigl(\xi^{u^\vee}(w^\vee)\bigr)$ where $\varphi$ was
defined in \S \ref{eqcohs}, and $\xi^u(w)$ is the coefficient
considered in \cite{Bi}, Thm. 4.}\end{footnote}.

For $\beta_{i_j} \in \Delta$, denote $s_{\beta_{i_j}}$ by
$s_{i_j}$ and let $s_{i_1}\cdot...\cdot s_{i_p}$ be a reduced word
decomposition for the representative $w^\vee$ in $W^P$. For each
$j$ between $1$ and $p$ define (following \cite{Bi} \S 4)
\begin{equation}\label{rj} r_w(j)=s_{i_1}\cdot ... \cdot
s_{i_{j-1}}(\beta_{i_j}) \end{equation}

\begin{prop}[\cite{Bi}, Lemma 4.1, \cite{Ku} Prop. 11.1.11]\label{cuwprop}
Let $u,w$ be two representatives in $W^P$
such that $ u \prec w$. Then the equivariant coefficient
$c_{u,w}^w$ is equal to
\begin{equation}\label{cuw} \sum r_w(j_1)\cdot ... \cdot r_w(j_k)
\end{equation} where the sum is over all ordered sequences ${j_1} < {j_2} < ... < {j_k} $
such that $s_{i_{j_1}}\cdot...\cdot s_{i_{j_k}}$ is a reduced word
decomposition for $u^\vee$.
\end{prop}

As a Corollary, we would like to note the following combinatorial
properties of the function $R(u,w;f)$, which are not at all clear
from its definition:
\begin{cor}\label{assign} Let $u,w$ be two representatives such that $u \prec
w$. Then \begin{enumerate}\item[(a)] $R(u,w;f)$ does not depend on
the choice of the assignment $f$, i.e. \[R(u,w;f)=R(u,w;f')\] for
any two assignments $f$ and $f'$. Denote this function by
$R(u,w)$.

\item[(b)] $R(u,w)$ is not equal to zero. Consequently, for any $u
\prec w$, the coefficient $c_{u,w}^w$ is not equal to zero.

\item[(c)] $R(u,w)/R(\overline{w}_0,w)=\sum r_w(j_1)\cdot ...
\cdot r_w(j_k)$ where the sum is as in the equation (\ref{cuw}).
\end{enumerate}
\end{cor}
\begin{proof} If $u=\overline{w}_0$, both (a) and (b) follow
from formulae (\ref{top}) and (\ref{uw}). For general $u$ such
that $u \prec w$, (a) follows from the fact that
$R(\overline{w}_0,w;f)$ does not depend on $f$, and that the
quotient $R(u,w;f')/R(\overline{w}_0,w;f)$ does not depend on $f$
and $f'$, being equal to the equivariant coefficient $c_{u,w}^w$.
For (b), by Lemma \ref{noz}, there exists an assignment $f$ such
that $R(u,w;f)$ is not equal to zero, hence this is true for all
assignments, by (a) just proved. Then, formula (\ref{uw}) implies
that $c_{u,w}^w$ must also be nonzero. (c) follows from Prop.
\ref{cuwprop}.
\end{proof}

{\it Remark:} 1. The fact that $c_{u,w}^w$ is not equal to zero if
$u \prec w$ can also be derived from Prop. \ref{cuwprop}, or using
the fact that $c_{u,w}^w$ is equal to the localization of the
Schubert class $\sigma(u)$ to the $T-$fixed point $wP \in G/P$
(cf. \cite{A}).

\section{Appendix - Proof of the Lemma \ref{F}}
The aim of this Appendix is to sketch the proofs of some of the
properties of the equivariant coefficients $c_{s(\beta),w}^w$,
defined in \S \ref{eqcohs}. These properties are needed in the
proof of the algorithm for the EQLR coefficients. We use the
notations of \S \ref{schubert}.

Recall from formula (\ref{d}) that if $s_{i_1}\cdot...\cdot
s_{i_k}$ (where $s_{i_j}=s_{\beta_{i_j}}$) is a reduced word
decomposition for $w \in W$, and $\beta=\beta_i$ is a simple root
in $\Delta$ then
\[\label{appd} D(s_\beta,w) = \sum_{i_j=i} s_{i_1}\cdot...\cdot
s_{i_{j-1}}(\beta) \]

We show first that $D(s_\beta,w) \geqslant 0$. For that, it is
enough to show that each term of its defining sum is nonnegative
(in the sense above). This follows immediately from the following
Lemma:
\begin{lemma}\label{sas} Let $w \in W$ and $\alpha$ a positive
root such that $l(ws_\alpha)>l(w)$. Then $w(\alpha)>0$.
\end{lemma}
\begin{proof} See the proposition in \cite{Hum}, \S
5.7.\end{proof} The second Lemma gives an equivalent definition of
$D(s_\beta,w)$. Recall that $\omega_\beta$ is the fundamental
weight corresponding to the positive simple root $\beta$.

\begin{lemma}\label{weight} The coefficient
$D(s_\beta,w)$ is equal to
\[ D(s_\beta,w)=\omega_\beta - w \omega_\beta \]
\end{lemma}
\begin{proof} See \cite{Ku},
Cor. 1.3.22. or Thm. 11.1.7(c). \end{proof}

\begin{prop}\label{neq} Let $u,w$ be two distinct representatives in
$W^P$. There exists a positive simple root $\beta$ in $\dnop$ such
that $D(s_\beta,w)-D(s_\beta,u)$ is not equal to zero.
\end{prop}
\begin{proof} Assume $D(s_\beta,w)=D(s_\beta,u)$ for any $\beta
\in \dnop$. Then Lemma \ref{weight} implies that
$w(\omega_\beta)=u(\omega_\beta)$, i.e.
$u^{-1}w(\omega_\beta)=\omega_\beta$ for all $\beta$ as before.
Take $\rho =\sum_{\beta \in \dnop} \omega_\beta$. Then
$(\rho,\beta)>0$ for any $\beta \in \dnop$ and $(\rho,\beta)=0$
for any $\beta \in \Delta_P$. Moreover, $u^{-1}w(\rho)=\rho$. Then
by \cite{B}, Ch. 5, \S 4.6 (see also \cite{Hum}, Thm. from \S
1.12) it follows that $u^{-1}w$ must be in $W_P$, which
contradicts the hypothesis.
\end{proof}

The next Proposition shows that the difference
$D(s_\beta,w)-D(s_\beta,u)$, for $\beta$ in $\Delta$, satisfies a
positivity property, provided that $u$ is less than $w$ in the
Bruhat ordering. Recall (see e.g. \cite{Hum}, \S 5.9) that in this
ordering, denoted $\leqslant$ , $u$ is less than $w$ if there is a
chain $u=u_0, u_1,...,u_k=w$ such that $u_{i+1}=u_is_\alpha$ where
$\alpha$ is a positive root such that $l(u_{i+1})>l(u_i)$. Note
that $u \leqslant w$ in the Bruhat ordering if and only if $w
\prec u$ in the ordering defined in \S \ref{formulaes}.
\begin{prop}\label{pos} Let $u,w$ be two permutations in $W$ such
that $u \leqslant w$, and let $\beta$ be a positive simple root.
Then $D(s_\beta,w)-D(s_\beta,u)\geqslant 0$. \end{prop}
\begin{proof} First note that one can reduce to the case when $w$
covers $u$, hence $w=us_\alpha$ for $\alpha$ a positive root. Then
$w(\omega_\beta)=us_\alpha(\omega_\beta)$ is equal to
\[u(\omega_\beta-(h_\alpha,\omega_\beta)\alpha)=u(\omega_\beta)-(h_\alpha,\omega_\beta)u(\alpha)\]
hence \begin{equation}\label{Ddif} D(s_\beta,w)-D(s_\beta,u) =
(h_\alpha,\omega_\beta)u(\alpha) \end{equation} Then
$(h_\alpha,\omega_\beta)\geqslant 0$ since $\alpha$ is a positive
root and $u(\alpha)>0$ by Lemma \ref{sas}, since
$l(us_\alpha)>l(u)$.\end{proof}

\begin{prop}\label{posP} Let $u,w$ be two distinct representatives in $W^P$
such that $u \leqslant w$ in Bruhat ordering. Let
$Cov_{\leqslant}(u,w)$ be the set of positive roots $\alpha$ in
$\pnop$ such that $us_\alpha$ is a representative in $W^P$,
$l(us_\alpha)=l(u)+1$ and $us_\alpha \leqslant w$. Then

 (a) The set $Cov_{\leqslant}(u,w)$ is nonempty.

(b) Let $\alpha$ be any positive root in $Cov_{\leqslant}(u,w)$
and $\beta$ a positive simple root in $\dnop$ such that
$D(s_\beta,us_\alpha)-D(s_\beta,u)$ is not equal to zero (such a
$\beta$ exists by Prop. \ref{neq}). Then
$h_\alpha(\omega_\beta)>0$.
\end{prop}
\begin{proof} By Prop. from  \S 5.11 in \cite{Hum} there exists a chain
\[ u=u_0 \leqslant u_1 \leqslant ...\leqslant u_k=w \] in $W$,
with $k=l(w)-l(u)$ such that $u_{i+1}=u_is_{\alpha_i}$ and
$l(u_{i+1})=l(u_i)+1$, with $\alpha_i$ in $\Phi^+$. Modulo the
Weyl group $W_P$ of $P$ this determines a chain in $W/W_P$ between
$uW_P$ and $wW_P$, necessarily of the same length $k$, since $u,w$
are representatives in $W^P$ whose length difference is $k$. In
particular, this shows that no $\alpha_i$ can be in $\Phi_P^+$ and
that the $u_i$'s are minimal length representatives for the cosets
in $W/W_P$, i.e. that $u_i$'s are in $W^P$. In particular, the
positive root $\alpha_0$ defined by $u_1=u_0s_{\alpha_0}$ must be
in $Cov_{\leqslant}(u,w)$, which finishes the proof of (a).

To prove (b), note that Prop. \ref{pos} implies that
$D(s_\beta,us_\alpha)-D(s_\beta,u)>0$, hence \[
h_\alpha(\omega_\beta) u(\alpha) > 0 \] by formula (\ref{Ddif}).
Lemma \ref{sas} implies that $u(\alpha)>0$ (because
$l(us_\alpha)>l(u)$), so $h_\alpha(\omega_\beta)$ must be a
nonzero, hence positive, integer, which ends the proof.
\end{proof}

We interpret now these properties in the terms of the equivariant
coefficients $c_{s(\beta),w}^w$ defined in Prop. \ref{eqc}. Recall
that these are defined as $\varphi(D(s_\beta,w^\vee))$ where
$\varphi$ is the automorphism of $\Phi$ sending the positive
simple root $\beta_i$ to the negative simple root $w_0(\beta_i)$
(to see this note that $w_0(\beta_i)$ must be a negative root, by
Lemma \ref{sas}, and that $s_{w_0(\beta_i)}$ is equal to
$w_0s_{\beta_i}w_0$; the last permutation sends all the positive
roots but one to positive roots, hence $l(w_0s_{\beta_i}w_0)=1$).
To state these properties, we recall from \S \ref{formulaes} the
analogue of the set $Cov_{\leqslant}(u,w)$ for the ordering
$\prec$. This set, denoted for simplicity $Cov(u,w)$ is the set of
positive roots $\alpha$ in $\pnop$ such that $us_\alpha$ is a
representative in $W^P$, $c(us_\alpha)=c(u)+1$ and $us_\alpha
\prec w$.

\begin{cor}[also Lemma \ref{F}]\label{properties} Let $u,w$ be two representatives in $W^P$.
\begin{enumerate}\item The coefficient $c_{s(\beta),w}^w$, for $\beta$ in $\dnop$,
is a linear homogeneous combination of negative simple roots with
nonnegative coefficients, and there exists a $\beta$ for which
this coefficient is nonzero.

\item There exists a positive simple root $\beta$ in $\dnop$ such
that $c_{s(\beta),w}^w-c_{s(\beta),u}^u$ is nonzero.

\item Assume that $u \prec w$ in the reverse Bruhat ordering
defined in \S \ref{formulaes}. Then for any $\beta$ in $\dnop$ the
difference $c_{s(\beta),w}^w-c_{s(\beta),u}^u$ is a nonnegative
combination of negative simple roots.

\item If $u \prec w$ the set $Cov(u,w) $ is nonempty. Moreover,
for any positive root $\alpha$ in $Cov(u,w)$ and any $\beta$
positive simple root in $\dnop$ such that
$c_{s(\beta),us_\alpha}^{us_\alpha}-c_{s(\beta),u}^u$ is nonzero,
the integer $h_\alpha(\omega_\beta)$ is positive.
\end{enumerate}
\end{cor}
\begin{proof} (1) follows from the similar property of
$D_{s_\beta,w}$, (2) from Prop. \ref{neq}, while (3) follows from
Prop. \ref{pos} taking into account that $u \prec w$ in the
reversed Bruhat ordering if and only if $u^\vee \leqslant w^\vee$
in the usual Bruhat ordering. Assertion (4) follows then from
Prop. \ref{posP}.
\end{proof}

\end{document}